\documentclass[a4paper]{article}
\usepackage{url}
\usepackage{amsfonts,amsmath}

\begin{document}

%
\title{Performance Evaluation of an Extrapolation Method for Ordinary Differential Equations with Error-free Transformation}
\author{Tomonori Kouya\\\url{kouya.tomonori@sist.ac.jp}\\Shizuoka Institute of Science and Technology\\2200-2 Toyosawa, Fukuroi, Shizuoka 437-8555 JAPAN}
\date{\today}
\maketitle
\abstract{
The application of error-free transformation (EFT) is recently being developed to solve ill-conditioned problems. It can reduce the number of arithmetic operations required, compared with multiple precision arithmetic, and also be applied by using functions supported by a well-tuned BLAS library. In this paper, we propose the application of EFT to explicit extrapolation methods to solve initial value problems of ordinary differential equations. Consequently, our implemented routines can be effective for large-sized linear ODE and small-sized nonlinear ODE, especially in the case when harmonic sequence is used.
}

\section{Introduction}
Double- or multi-fold arithmetic, which is implemented using error-free transformation\cite{eft} (EFT), is recently being paid attention to be parallel with multiple precision arithmetic. It can reduce the number of normalizations occurring in each multiple precision arithmetic and also be applied using functions supported by well-tuned BLAS libraries such as the Intel Math Kernel library and OpenBLAS. Kobayashi and Ogita\cite{kobayashi2015} have demonstrated the effectiveness of double-fold arithmetic using matrix arithmetic provided by BLAS level 3 (BLAS3) to solve ill-conditioned linear equations.

In this paper, we propose the application of EFT to explicit extrapolation methods to solve initial value problems (IVPs) of ordinary differential equations (ODEs). The explicit extrapolation methods can be implemented only using vector arithmetic provided by BLAS level 1 (BLAS1). We implemented double-fold explicit extrapolation methods and evaluated the performance of various precision techniques such as double and double-double (DD) arithmetic and algorithms such as classical M\o ller method to reduce the accumulation of round-off errors. Consequently, our implemented routines can be effective for large-sized linear ODE and small-sized nonlinear ODE, especially when harmonic sequence is used.

\section{Explicit Extrapolation for ODEs}
The $n$-dimensional IVP and ODE to be solved is shown as follows:
\begin{gather}
\left\{\begin{split}
	\frac{d\mathbf{y}}{dt} &= \mathbf{f}(t, \mathbf{y}) \\
	\mathbf{y}(t_\mathrm{start}) &= \mathbf{y}_\mathrm{start}
\end{split}\right. \label{eqn:ode_ivp} \\
\mbox{Integration interval:}\ [t_\mathrm{start}, t_\mathrm{end}] \ni t, \nonumber
\end{gather}
where $\mathbf{y}, \mathbf{f}(t, \mathbf{y}) \in \mathbb{R}^n$. We discretize the above integration interval, and at each $t_{\rm next}$ $\in [t_\mathrm{start}, t_{\rm end}]$, compute the approximation $\mathbf{y}_{\rm next} \approx \mathbf{y}(t_{\rm next})$ from $\mathbf{y}_{\rm old} \approx \mathbf{y}(t_{\rm old})$ using the explicit extrapolation method. In the rest of this section, we describe the algorithm in detail and the propagation of round-off errors in the extrapolation process shown in the Hairer \& Wanner's textbook\cite{hairer}.

\subsection{Algorithm of explicit extrapolation method}
At first, we set maximum number of stages as $L$, relative tolerance $\varepsilon_R$, absolute tolerance $\varepsilon_A$ and support sequence $\{w_i\}_{i=1}^L$. We use two types of support sequence: Romberg sequence ($w_i := 2^i$) or harmonic sequence ($w_i := 2(i + 1)) $.

The standard explicit extrapolation method uses a combination of the explicit Euler method
\begin{equation}
	\mathbf{y}_{1} := \mathbf{y}_0 + h \mathbf{f}(t_0, \mathbf{y}_0) = \mathbf{y}_0 + h \mathbf{f}_0, \label{eqn:explicit_euler}
\end{equation}
and the mid-point methods
\begin{equation}
\begin{split}
		\mathbf{y}_{k + 1} &:= \mathbf{y}_{k-1} + 2 h \mathbf{f}(t_k, \mathbf{y}_k) \\
		&= \mathbf{y}_{k-1} + 2 h \mathbf{f}_k\ (k = 1, 2,..., w_i - 1).
	\end{split} \label{eqn:explicit_mid_point}
\end{equation}
Through the above process, we can obtain the initial sequence as $\mathbf{T}_{i1} := \mathbf{y}_{w_i}$. The step size $h$ is determined as $ h := (t_{\rm next} - t_{\rm old}) / w_i$. Then, each discretization point $t_k$ is fixed as $ t_k := t_{\rm old} + k h \in [t_{\rm old}, t_{\rm next}]$ in the above process to obtain the initial sequence, where $t_0 := t_{\rm old},\ \mathbf{y}_0 \approx \mathbf{y}(t_{\rm old})$.

Next, we calculate $\mathbf{T}_{ij}$$(j = 2,..., i)$ using $\mathbf{T}_{i-1,j-1}$ and $\mathbf{T}_{i,j-1}$ as follows:
\begin{gather}
	c_{ij}  := \left(\left(\frac{w_{i}}{w_{i-j+1}}\right)^2 - 1\right)^{-1} \nonumber \\
	\mathbf{R}_{ij} :=  c_{ij}(\mathbf{T}_{i, j-1} - \mathbf{T}_{i-1,j-1})\ \label{eqn:extrap_process} \\
	\mathbf{T}_{ij} :=  \mathbf{T}_{i,j-1} + \mathbf{R}_{ij}. \nonumber
\end{gather}

In the above extrapolation process, we must check if the following convergence condition is satisfied.
\begin{equation}\begin{split}
	\|\mathbf{R}_{ij}\| &\le \varepsilon_R\|\mathbf{T}_{i,j-1}\| + \varepsilon_A \\
\end{split} \label{eqn:extrap_converge}
\end{equation}
If satisfied, we fix $\mathbf{y}_{\rm next} := \mathbf{T}_{ij}$; if not, we calculate additional approximation $\mathbf{T}_{i+1,1}$ from the process of (\ref{eqn:explicit_euler}) and (\ref{eqn:explicit_mid_point}) and continue the extrapolation process (\ref{eqn:extrap_process}). This iteration is definitely stopped at $i, j = L$ if not converged. Murofushi and Nagasaka\cite{murofushi_nagasaka1994} proposed $\varepsilon_R = \varepsilon_A = 0$ as tolerances in (\ref{eqn:extrap_converge}) to obtain the optimized approximation when global truncation and round-off error are balanced.

This above-mentioned explicit extrapolation method can be implemented by using AXPY and SCAL(\figurename\ \ref{fig:blas1}) supported in BLAS1.

%
\begin{figure}[htbp]
\begin{center}
\begin{minipage}{.35\textwidth}
\begin{tabular}{l}\hline
	$\mathbf{y} := \mbox{AXPY}(\alpha, \mathbf{x}, \mathbf{y})$ \\ \hline
	\hspace{2mm} $\mathbf{y} := \alpha\otimes \mathbf{x} \oplus \mathbf{y}$ \\
	\hspace{2mm} return $\mathbf{y}$ \\ \hline
\end{tabular}
\end{minipage}
\begin{minipage}{.35\textwidth}
\begin{tabular}{l}\hline
	$\mathbf{x} := \mbox{SCAL}(\alpha, \mathbf{x})$ \\ \hline
	\hspace{2mm} $\mathbf{x} := \alpha\otimes \mathbf{x}$ \\
	\hspace{2mm} return $\mathbf{x}$ \\ \hline
\end{tabular}
\end{minipage}
\end{center}
\caption{Standard BLAS1 functions: AXPY and SCAL}\label{fig:blas1}
\end{figure}

\subsection{Propagation of Round-off error in extrapolation process}

Hairer \& Wanner\cite{hairer} analyzed the propagation effect of round-off error in the extrapolation process (\ref{eqn:extrap_process}). The following assumptions are supposed.
\begin{itemize}
\item The initial sequence $T_{i1}$ contains $\varepsilon_{i1} = (-1)^{i-1}\varepsilon$ as the corresponding error.
\item These errors do not diminish each other in the extrapolation process (\ref{eqn:extrap_process}).
\end{itemize}
According to (\ref{eqn:extrap_process}), the error $\varepsilon_{ij}$ in $\mathbf{T}_{ij}$ is expressed as
\begin{equation}
	\varepsilon_{ij} = \varepsilon_{i,j-1} + c_{ij} (\varepsilon_{i, j-1} - \varepsilon_{i-1,j-1}) = r_{ij}\varepsilon. \label{eqn:roundoff_error_propagation}
\end{equation}
The coefficient $r_{ij}$ in $\varepsilon_{ij} = r_{ij}\varepsilon$ uncovers the propagation effect of round-off errors in the initial sequence. In case of $L=20$, the effect is less than two times using Romberg sequence. On the contrary, the extrapolation process may provide a $O(10^6)$ propagation effect using harmonic sequence.

Murofushi and Nagasaka\cite{murofushi_nagasaka1994} recommend choosing Romberg sequence as the support sequence to limit the propagation effect of round-off errors. Although harmonic sequence may increases errors in approximation, it can reduce the number of calculations required for obtaining the initial sequence. Therefore, it can get better performance with harmonic sequence as compared with Romberg sequence when heavy multiple precision arithmetic is applied.

\section{Explicit extrapolation method with error-free transformation}

As already described, the explicit extrapolation method can be implemented only by using SCAL and AXPY of BLAS1 functions. In this section, we extend the two functions to these ones with error evaluations using EFT and then describe the explicit extrapolation method with error evaluation by these extended BLAS1 functions.

\subsection{BLAS1 functions with EFT}

We denote standard IEEE754 elementary arithmetic operators as $\oplus$, $\otimes$, $\ominus$, and $\oslash$. For these arithmetic operators, we call error-free transformation (EFT), which these functions in \figurename\ \ref{fig:basic_function_eft} can provide the corresponding errors occurring in these elementary arithmetic operators.

%
\begin{figure}[htbp]
\begin{center}
\begin{tabular}{l}\hline
	($s$, $e$) $:=$ QuickTwoSum($a$, $b$) \\ \hline
	\hspace{2mm} $s := a \oplus b$;\ $e := b \ominus (s \ominus a)$ \\
	\hspace{2mm} return ($s$, $e$) \\ \hline
\end{tabular}

\begin{tabular}{l}\hline
	($s$, $e$) $:=$ TwoSum($a$, $b$) \\ \hline
	\hspace{2mm} $s := a \oplus b$;\ $v := s \ominus a$ \\
	\hspace{2mm} $e := (a \ominus (s \ominus v)) \oplus (b \ominus v)$ \\
	\hspace{2mm} return ($s$, $e$) \\ \hline
\end{tabular}

\begin{tabular}{l}\hline
	($s$, $e$) $:=$ TwoProd($a$, $b$) (with FMA) \\ \hline
	\hspace{2mm} $s := a \otimes b$ \\
	\hspace{2mm} $e := \mbox{FMA}(a, b, -s)$\ ($= a \times b - s$) \\
	\hspace{2mm} return ($s$, $e$) \\ \hline
\end{tabular}

\caption{Basic functions of Error-Free Transformation}\label{fig:basic_function_eft}
\end{center}
\end{figure}

For implementing of SCAL and AXPY with error evaluations, FMA arithmetic with errors is desirable. We use the FMAerror function (\figurename\ \ref{fig:fma_eft}) proposed by S.Boldo \& J-M. Muller\cite{fma_error}.

%
\begin{figure}[htbp]
\begin{center}

\begin{tabular}{l}\hline
	($s$, $e_1$, $e_2$) $:=$ FMAerror($a$, $x$, $y$) \\ \hline
	\hspace{2mm} $s := \mbox{FMA}(a, x, y)$;\	$(u_1, u_2) := \mbox{TwoProd}(a, x)$ \\
	\hspace{2mm} $(\alpha_1, \alpha_2) := \mbox{TwoSum}(y, u_2)$;\ $(\beta_1, \beta_2) := \mbox{TwoSum}(u_1, \alpha_1)$ \\
	\hspace{2mm} $\gamma  := \beta_1 \ominus s \oplus \beta_2$;\ $(e_1, e_2) := \mbox{QuickTwoSum}(\gamma, \alpha_2)$ \\
	\hspace{2mm} return ($s$, $e_1$, $e_2$) \\ \hline
\end{tabular}
\caption{FMA arithmetic with error evaluation}\label{fig:fma_eft}
\end{center}

\end{figure}

FMAerror guarantees $s + e_1 + e_2 = ax + y$, where $s = a\otimes x \oplus y$, $|e_1 + e_2| = (1/2)\mathbf{u}|s|$ ($\mathbf{u}$ is unit of round-off error), and $|e_2| = \frac{1}{2} \mathbf{u} |e_1|$.

A similar function with FMAerror can be implemented using Sloppy DDadd and DDmul operators that are supported in DD libraries. As shown in \tablename\ \ref{table:compare_fma}, the total number of elementary arithmetic operations is the same. In our implementation, we use FMAerror to implement our BLAS1 functions with error evaluations.

%
\begin{table}[htbp]
\begin{center}
\caption{Number of elementary arithmetic}\label{table:compare_fma}
\begin{tabular}{c|c|c|c}
										& $\oplus$, $\ominus$ & $\otimes$ & FMA \\ \hline
FMAerror						&    17         &    1      &  2  \\ \hline
Sloppy DDadd \& DDmul	&    16         &    3      &  1  \\ \hline
\end{tabular}
\end{center}
\end{table}

Using basic EFT arithmetic, AXPYerror and SCALerror can be implemented as shown in \figurename\ \ref{fig:blas1_error}, where each EFT arithmetic is applied for each element of the vectors.

%
\begin{figure}[htbp]
\begin{center}

\begin{tabular}{l}\hline
	$(\mathbf{y}, \mathbf{e}_\mathbf{y}) := \mbox{AXPYerror}(\alpha, e_\alpha, \mathbf{x}, \mathbf{e}_\mathbf{x}, \mathbf{y}, \mathbf{e}_\mathbf{y})$ \\ \hline
	\hspace{2mm} $(\mathbf{y}, \mathbf{e}_1, \mathbf{e}_2) := \mbox{FMAerror}(\alpha, \mathbf{x}, \mathbf{y})$ \\
	\hspace{2mm} $\mathbf{e}_\mathbf{y} := \mathbf{e}_1 \oplus \mathbf{e}_2 \oplus \alpha \otimes \mathbf{e}_\mathbf{x} \oplus e_\alpha \otimes \mathbf{x} \oplus \mathbf{e}_\mathbf{y}$ \\
	\hspace{2mm} return ($\mathbf{y}$, $\mathbf{e}_\mathbf{y}$) \\ \hline
\end{tabular}

\begin{tabular}{l}\hline
	$(\mathbf{x}, \mathbf{e}_\mathbf{x}) := \mbox{SCALerror}(\alpha, e_\alpha, \mathbf{x}, \mathbf{e}_\mathbf{x})$ \\ \hline
	\hspace{2mm} $(\mathbf{w}_1, \mathbf{w}_2) := \mbox{TwoProd}(\alpha, \mathbf{x})$ \\
	\hspace{2mm} $\mathbf{w}_2 := \alpha \otimes \mathbf{e}_\mathbf{x} \oplus e_\alpha \otimes (\mathbf{x} \oplus \mathbf{e}_\mathbf{x}) \oplus \mathbf{w}_2 $ \\
	\hspace{2mm} $(\mathbf{x}, \mathbf{e}_\mathbf{x}) := \mbox{QuickTwoSum}(\mathbf{w}_1, \mathbf{w}_2)$ \\
	\hspace{2mm} return ($\mathbf{x}$, $\mathbf{e}_\mathbf{x}$) \\ \hline
\end{tabular}
\caption{BLAS1 with error evaluation: AXPYerror and SCALerror}\label{fig:blas1_error}
\end{center}
\end{figure}

%
\subsection{Explicit extrapolation method with EFT}

We can implement the explicit extrapolation method with EFT using BLAS1 functions with error evaluations. Suppose that we can evaluate $\mathbf{f}(t_k + e_{t_k}, \mathbf{y}_k + \mathbf{e}_\mathbf{y_k}) = \mathbf{f}_k + \mathbf{e}_{\mathbf{f}_k}$ with its error.

The explicit Euler method (\ref{eqn:explicit_euler}) is extended as
\begin{equation}
\begin{split}
		(\mathbf{y}_1, \mathbf{e}_{\mathbf{y}_1}) &:= (\mathbf{y}_0, \mathbf{e}_{\mathbf{y}_0}) \\
		(\mathbf{y}_{1}, \mathbf{e}_{\mathbf{y}_1}) &:= \mbox{AXPYerror}(h, e_h, \mathbf{f}_0, \mathbf{e}_{\mathbf{f}_0}, \mathbf{y}_1, \mathbf{e}_{\mathbf{y}_1}).
\end{split} \label{eql:explicit_euler_eft}
\end{equation}

The explicit mid-point method (\ref{eqn:explicit_mid_point}) is extended as
\begin{equation}
\begin{split}
		(\mathbf{y}_{k+1}, \mathbf{e}_{\mathbf{y}_{k+1}}) &:= (\mathbf{y}_{k-1}, \mathbf{e}_{\mathbf{y}_{k-1}}) \\
		(\mathbf{y}_{k+1}, \mathbf{e}_{\mathbf{y}_{k+1}}) &:= \mbox{AXPYerror}(2\otimes h, 2\otimes e_h, \mathbf{f}_{k}, \\
		&\ \mathbf{e}_{\mathbf{f}_{k}}, \mathbf{y}_{k+1}, \mathbf{e}_{\mathbf{y}_{k+1}})\ (k = 1, 2, ..., w_i - 1).
\end{split} \label{eqn:explicit_mid_point_eft}
\end{equation}

Therefore, the initial sequence is obtained as $(\mathbf{T}_{i1}, \mathbf{e}_{\mathbf{T}_{i1}})$ $:= (\mathbf{y}_{w_{i}}, $$ \mathbf{e}_{\mathbf{y}_{w_i}})$.

For the preparation of the extrapolation process, we calculate $c_{ij}$ in (\ref{eqn:extrap_process}) as $(c_{ij}, e_{c_{ij}}) := 1 / ((w_{i} / w_{i-j+1})^2 - 1)$ by application of the DD arithmetic.

Extrapolation process (\ref{eqn:extrap_process}) is extended as
\begin{equation}
\begin{split}
	(\mathbf{T}_{ij}, \mathbf{e}_{\mathbf{T}_{ij}}) &:= (\mathbf{T}_{i, j-1}, \mathbf{e}_{\mathbf{T}_{i, j-1}}) \\
	(\mathbf{R}_{ij}, \mathbf{e}_{\mathbf{R}_{ij}}) &:= (\mathbf{T}_{i, j-1}, \mathbf{e}_{\mathbf{T}_{i, j-1}}) \\
	(\mathbf{R}_{ij}, \mathbf{e}_{\mathbf{R_{ij}}}) &:= \mbox{AXPYerror}(-1, 0, \mathbf{T}_{i-1, j-1}, \mathbf{e}_{\mathbf{T}_{i-1, j-1}}, \\
		&\ \mathbf{R}_{ij}, \mathbf{e}_{\mathbf{R}_{ij}}) \\
	(\mathbf{R}_{ij}, \mathbf{e}_{\mathbf{R}_{ij}}) &:= \mbox{SCALerror}(c_{ij}, e_{c_{ij}}, \mathbf{R}_{ij}, \mathbf{e}_{\mathbf{R}_{ij}}) \\
	(\mathbf{T}_{ij}, \mathbf{e}_{\mathbf{T}_{ij}}) &:= \mbox{AXPYerror}(1, 0, \mathbf{R}_{ij}, \mathbf{e}_{\mathbf{R}_{ij}}, \mathbf{T}_{ij}, \mathbf{e}_{\mathbf{T}_{ij}}).
\end{split} \label{eqn:extrap_process_eft}
\end{equation}

%
\subsection{M\o ller method}

The M\o ller method is proposed to reduce accumulation of round-off errors incurred during approximation of IVPs of ODEs and is a type of compensated summation. For the original summation $S_i := S_{i-1} + z_{i-1}$, we compute it as follows:
\begin{eqnarray*}
	s_i &:=& z_{i-1} \ominus R_{i-1}\ (R_0 = 0) \\
	S_i &:=& S_{i-1} \oplus s_i;\ r_i := S_i \ominus S_{i-1};\ R_i := r_i \ominus s_i.
\end{eqnarray*}
The above formula can be rewritten using $R_{i}' = -R_i$ and QuickTwoSum in \figurename\ \ref{fig:basic_function_eft} as follows:
\begin{equation}
\begin{split}
	s_i &:= z_{i-1} \oplus R_{i-1}'\ (R_0' = 0) \\
	(S_i, R_i') &:= \mbox{QuickTwoSum}(S_{i-1}, s_i). \\
\end{split} \label{eqn:moeller_quick_two_sum}
\end{equation}
QuickTwoSum($S_{i-1}, s_i$) can obtain the correct error only in the case of $|S_{i-1}| \geq |s_i|$. Such situations can be expected in the process obtaining the initial sequence and in the extrapolation process, when the effect of the round-off error is larger than the  truncation error. Although the situation of being able to provide the correct error may be satisfied in practical situations, the effectiveness of the application of the M\o ller method is not observed in some cases. For our comparison, we apply the formula (\ref{eqn:moeller_quick_two_sum}) as M\o ller methods to (\ref{eqn:explicit_euler}), (\ref{eqn:explicit_mid_point}), and (\ref{eqn:extrap_process}).

\section{Numerical Experiments}
%
We compare the performances and relative errors using our implemented explicit extrapolation methods. Our computational environment is as follows:
\begin{description}
\item[H/W] AMD Ryzen 1700 (2.7 GHz), 32 GB RAM
\item[S/W] Ubuntu 16.04.5 x86\_64, GCC 5.4.0, QD 2.3.18\cite{qd}, LAPACK 3.8.0.
\end{description}

Our targets of precision are IEEE754 double precision (Double) and DD provided by the QD library, and the targeted algorithms are as follows:
\begin{description}
	\item[DEFT]:\ Double precision\ (\ref{eql:explicit_euler_eft}), (\ref{eqn:explicit_mid_point_eft}), (\ref{eqn:extrap_process_eft}), and $\mathbf{f} + \mathbf{e}_{\mathbf{f}}$
	\item[DEFT2]:\ Double precision\ (\ref{eql:explicit_euler_eft}), (\ref{eqn:explicit_mid_point_eft}), (\ref{eqn:extrap_process_eft}), $\mathbf{f}$, $\mathbf{e}_{\mathbf{f}} := 0$
	\item[DM\o ller]:\ Double precision M\o ller method.
\end{description}

DEFT2 means usage of the double precision $\mathbf{f}$ namely the error term of $\mathbf{f}$ is zero. For DEFT and DD computations, we use DD precision $\mathbf{f}$. For checking convergence (\ref{eqn:extrap_converge}), we use $\varepsilon_R = \varepsilon_A = 0$ unless otherwise specified.  All EFT basic functions are coded as C macros.

%
\subsection{Homogeneous linear ODE}

We pick up 2048-dimensional homogeneous linear ODE only using BLAS1 functions as follows:

\[ \begin{split}
\displaystyle\frac{d\mathbf{y}}{dt} &= [ -y_1\ \cdots\ -ny_n]^T \\
\mathbf{y}(0) &= [1\ \cdots\ 1]^T,\ t \in [0, 1/4].
\end{split} \]

The analytical solution is $\mathbf{y}(t)$ $ = [ \exp(-t)\ $ $\cdots\ $ $\exp(-nt)]^T$. This is simply one, so we use fixed step sizes $t_{\rm next} - t_{\rm old}$ $:= (1/4) $ $/ (\mbox{\#steps})$ for all patterns. 

%

\tablename\ \ref{table:romberg_l4_lode} shows the computational time (Unit: s) and its maximum relative errors for all elements of approximation at $t_{\rm end} = 1/4$ in the case of Romberg sequence and $L=4$. The line of the table shows the boundary being at the same level as that of maximum relative errors.

\begin{table}
\begin{center}
\caption{Linear ODE: Romberg sequence: $L=4$ at $t_{\rm end} = 1/4$}\label{table:romberg_l4_lode}
\begin{tabular}{|c|ccccc|ccccc|}\hline
$L=4$	& \multicolumn{5}{|c|}{Computational time (s)} \\
\#steps & DD	& DEFT&	DEFT2	&Double	& DM\o ller	\\ \hline
512		&1.79	&1.41	&1.4	&0.2	&0.33		\\
1024	&3.59	&2.81	&2.82	&0.41	&0.67		\\ \cline{2-6}
2048	&7.18	&5.64	&5.64	&0.81	&1.33		\\
4096	&14.4	&11.3	&11.3	&1.62	&2.66		\\ 
8192	&28.8	&22		&22		&3.17	&5.33		\\ \hline
\#steps & \multicolumn{5}{|c|}{Max. Relative Error} \\
512		&1.84E-07	&1.8E-07	&1.8E-07	&1.8E-07	&1.8E-07\\
1024	&1.17E-10	&1.2E-10	&1.2E-10	&1.2E-10	&1.2E-10\\ \cline{2-6}
2048	&9.28E-14	&9.3E-14	&9.4E-14	&1.5E-13	&9.4E-14\\
4096	&8.18E-17	&4.6E-16	&1.6E-14	&2.3E-13	&4.3E-14\\ 
8192	&7.59E-20	&3.3E-16	&2.4E-14	&3.9E-13	&1.7E-13\\ \hline
\end{tabular}
\end{center}
\end{table}

Consequently, we can observe the following results.
\begin{itemize}
	\item At the same order of maximum relative error, DEFT is approximately 1.3 times faster than DD. The difference in performance between DEFT and DEFT2 cannot be observed.
	\item Except DD over 2048 \#steps, the relative error of DEFT is the smallest. The M\o ller method reduced the relative error by approximately 1 decimal digit, when compared with double precision.
\end{itemize}

%
\tablename\ \ref{table:harmonic_l6_lode} shows the computational time and its maximum relative errors in the case of harmonic sequence and $L=6$.
\begin{table}
\begin{center}
\caption{Linear ODE: Harmonic sequence: $L=6$ at $t_{\rm end} = 1/4$}\label{table:harmonic_l6_lode}
\begin{tabular}{|c|ccccc|ccccc|}\hline
$L=6$	& \multicolumn{5}{|c|}{Computational Time (s)} \\
\#steps & DD	& DEFT&	DEFT2	&Double	& DM\o ller	\\ \hline
512		&1.87	&1.76	&1.31	&0.28	&0.4	\\  \cline{2-6}
1024	&3.74	&3.53	&2.63	&0.55	&0.81	\\
2048	&7.48	&6.93	&5.25	&1.11	&1.62	\\
4096	&14.9	&10.4	&10.5	&2.22	&3.24	\\ 
8192	&29.9	&15.4	&21	&4.43	&6.49	\\ \hline
\#steps& \multicolumn{5}{|c|}{Max. Relative Error} \\
512		&4.3E-10	&4.3E-10	&4.3E-10	&4.3E-10	&4.3E-10\\ \cline{2-6}
1024	&1.7E-14	&2.7E-14	&2.7E-14	&7.1E-13	&6.6E-13\\
2048	&8.4E-19	&1.3E-14	&1.4E-14	&9.2E-13	&7.2E-13\\
4096	&4.6E-23	&5.5E-15	&1.1E-14	&1.0E-12	&7.6E-13\\ 
8192	&2.7E-27	&2.2E-15	&7.4E-15	&1.5E-12	&8.6E-13\\ \hline
\end{tabular}
\end{center}
\end{table}

Consequently, we can observe the following results.
\begin{itemize}
	\item At the same order of relative errros, DEFT's performance is approximately 6\%--7\% better than DD and is the same as DEFT2. Faster convergence than DEFT2 can be observed at \#steps$ = 8192$. 
	\item The smallest relative error can be obtained by DEFT without DD.
\end{itemize}

The above numerical experiments for homogeneous linear ODE demonstrate that DEFT performs better than DD for the same level of relative errors. DEFT2 can obtain better approximation than DM\o ller but cannot get better performance than DEFT.

%
\subsection{Resonance problem}

We pick up the following resonance problem that is necessary to control step sizes.

\[\begin{split}
\frac{d}{dt}\left[\begin{array}{c}
	y_1 \\
	y_2
\end{array}\right] &= \left[\begin{array}{c}
	y_2 \\
	-\alpha y_1^2 \sin t + 2\alpha y_1 y_2 \cos t
\end{array}\right] \\
\mathbf{y}(0) &= [1\ \alpha]^T,\ t \in [0, 37]
\end{split}\]
where $\alpha = 0.99999999$. The analytical solution is
\[ \left[\begin{array}{c}
	y_1 \\
	y_2
\end{array}\right] = \left[\begin{array}{c}
	1/(1-\alpha \sin t) \\
	\alpha \cos t / (1 - \alpha \sin t)^2
\end{array}\right]. \]

The algorithm of step size control is the same one proposed in Murofushi and Nagasaka\cite{murofushi_nagasaka1994}, wherein the current step size is halved if the convergent condition (\ref{eqn:extrap_converge}) is not satisfied. The maximum stages are $L = 12$ for Romberg sequence and $L=18$ for harmonic sequence as recommended in \cite{murofushi_nagasaka1994}.

%
\begin{table}
\begin{center}
\caption{Resonance Problem: Computational time and maximum relative errors at $t_{\rm end} = 37$}\label{table:resonance_comp_time_relerr}
\begin{tabular}{|c|c|c|c|}\hline
Romberg, $L=12$ &	\#steps & Comp.Time (s)&Max.Rel.Err. \\ \hline
DD($\varepsilon_R=10^{-16}$)		& 100 & 0.19 &  3.6E-04 \\
DEFT		& 100	& 0.42	& 3.7E-04 \\
DEFT2		& 84	& 0.04	& 3.5E-02 \\
Double		& 84	& 0.02	& 1.0E-01 \\ 
DM\o ller	& 98	& 0.07	& 5.2E-04 \\ \hline
\end{tabular}
\begin{tabular}{|c|c|c|c|}\hline
Harmonic, $L=18$	& \#steps & Comp.Time (s)	& Max.Rel.Err.\\ \hline
DD($\varepsilon_R=10^{-18}$)		& 186 & 0.06 &  6.0E-05 \\
DEFT		& 159	& 0.05	& 4.5E-04 \\ \hline
\end{tabular}
\end{center}
\end{table}

For any cases with Romberg sequence, we can obtain approximations at $t_{\rm end} = 37$ without breakdown; then, DEFT, DM\o ller. and DD($\varepsilon_R=10^{-16}$, $\varepsilon_A = 0$) can obtain the most precise approximations. DM\o ller's performance is the best in the case of Romberg sequence.

On the contrary, DD ($\varepsilon_R=10^{-18}$, $\varepsilon_A = 0$) and DEFT can obtain approximations without breakdown in the case of harmonic sequence. The DEFT with harmonic sequence can demonstrate the best performance at the same order of maximum relative error through all precision arithmetic and algorithms.

\section{Conclusion and future work}
We can conclude that the explicit extrapolation method with EFT is competitive for the DD arithmetic one. In future studies, we will implement and evaluate implicit extrapolation methods with EFT and its variation with BLAS2 and BLAS3 functions in various computational environments.

\end{document}